\documentclass[letterpaper, 10 pt, conference]{ieeeconf}  % Comment this line out if you need a4paper

\IEEEoverridecommandlockouts                              % This command is only needed if 
\usepackage{cite}
\usepackage{amsmath,amssymb,amsfonts}
\usepackage{algorithm}
\usepackage{algpseudocode}
\usepackage{graphicx}
\usepackage{textcomp}
\usepackage{caption}
\usepackage{color}
\usepackage{verbatim}

\usepackage{caption}
\newcommand{\be}{\begin{equation}}
\newcommand{\ee}{\end{equation}}
\newcommand{\bpsi}{{\boldsymbol \psi}}
\newcommand{\bPsi}{{\boldsymbol \Psi}}

\newcommand{\bw}{{\boldsymbol w}}
\newcommand{\bu}{{\boldsymbol u}}
\newcommand{\bz}{{\boldsymbol z}}

\newcommand{\lba}{\left[ \begin{array}}
\newcommand{\ear}{\end{array} \right]}
\DeclareMathOperator*{\argmin}{arg\,min}

\newtheorem{proposition}{Proposition}
\newtheorem{definition}{Definition}
\newtheorem{lemma}{Lemma}

\newtheorem{theorem}{Theorem}

\newtheorem{remark}{Remark}

\title{\LARGE \bf
A System Level Approach to LQR Control of the Diffusion Equation 
}

\author{Addie McCurdy, Andrew Gusty, and Emily Jensen% <-this % stops a space
\thanks{*This work was not supported by any organization}% <-this % stops a space
\thanks{A. McCurdy and A. Gusty are with the Department of Applied Mathematics, University of Colorado, Boulder.
        {\tt\small addie.mccurdy@colorado.edu, angu8719@colorado.edu}}%
\thanks{E. Jensen is with the Department of Electrical Computer \& Energy Engineering, University of Colorado, Boulder.
        {\tt\small ejensen@colorado.edu}}%
}

\begin{document}

\maketitle
\thispagestyle{empty}
\pagestyle{empty}

%%%%%%%%%%%%%%%%%%%%%%%%%%%%%%%%%%%%%%%%%%%%%%%%%%%%%%%%%%%%%%%%%%%%%%%%%%%%%%%%
\begin{abstract}
The optimal controller design problem for a linear, first-order spatially-invariant distributed parameter system is considered. Through a case study of the Linear Quadratic Regulator (LQR) problem for the diffusion equation over the torus, it is illustrated that the optimal controller design problem can be equivalently formulated as an optimization problem over the system's closed-loop mappings, analogous to the System Level Synthesis framework. This reformulation is solved analytically to recover the LQR for the diffusion equation, and an internally stable implementation of this controller is recovered from the optimal closed-loop mappings. It is further demonstrated that a class of spatio-temporal constraints on the closed-loop maps can be imposed on this closed-loop formulation while preserving convexity.

%The optimal controller design problem for a linear distributed parameter system subject to spatio-temporal controller locality constraints is typically non convex. Through a case study of LQR design for the diffusion equation over the unit circle, this work examines an approach for convex restrictions of this spatio-temporally constrained controller design problem motivated by the closed-loop design procedure of the System Level Synthesis framework. It is demonstrated that the unconstrained LQR problem for the diffusion equation can be reformulated as an infinite-dimensional optimization problem over the closed-loop mappings, with cost functional and constraints each affine in these mappings. An analytic solution to this reformulation is obtained in the spatial frequency domain and is shown to recover the well-known LQR solution for the diffusion equation. A non-minimal implementation of the resulting controller that preserves the spatio-temporal properties of the closed-loop maps is proposed, and this implementation is rigorously shown to provide internal stability. It is further demonstrated that spatio-temporal constraints on the closed-loop maps can be imposed in the reformulation of the controller design problem while preserving convexity.
\end{abstract}

%%%%%%%%%%%%%%%%%%%%%%%%%%%%%%%%%%%%%%%%%%%%%%%%%%%%%%%%%%%%%%%%%%%%%%%%%%%%%%%%

\section{Introduction}

As large-scale spatially-distributed systems with immense sensing and actuation become more prevalent, so does the need for tractable methods for optimal distributed controller design. One motivating example is the development of RF metamaterials for quantum computing, beamforming, and secure communications \cite{yang2023integrated, jung2024recent, wang2024reconfigurable}, which are  modeled as continnum structures governed by PDEs. To leverage distributed sensing and actuation in these and other PDE systems, it is necessary to design control policies that use minimal communication between distributed subcontroller units, each with access to a limited subset of system information. This requires solving an optimal control problem subject to controller locality constraints, which is known to be non-convex except in special settings \cite{rotkowitz2005characterization,bamieh2005convex}. 

To overcome the computational complexity of such constrained optimization problems, various convex restrictions or relaxations have been developed, e.g., \cite{sabuau2023network, zheng2020equivalence}. Of particular significance is the 
 \emph{System Level Synthesis} (SLS) framework \cite{wang2019system}, which is characterized by (1) a synthesis step to formulate the design problem and constraints as affine in the system's closed-loop mappings; and (2) an implementation step to construct an internally stable implementation of the controller from an interconnection of the optimal closed-loop maps. This implementation ensures that the controller inherits any structural properties constrained on the closed-loop maps. For example, imposing that the closed-loops are spatially sparse ensures that the resulting controller can also be implemented sparse elements.

The SLS framework was first formulated for finite-dimensional discrete-time systems over a finite time horizon \cite{wang2019system}, but has since been adapted to various other settings including infinite time-horizons \cite[Sec. 4]{anderson2019system}, countably infinite spatial dimensions \cite{jensen2018optimal}, and continuous-time dynamics \cite{jensen2021explicit}. However, these approaches could not be directly applied to PDE systems formulated over a continuous spatial domain, as noted by the recent work of Conger et. al. \cite{conger2025convex}. Indeed, it is well-known that discretized approaches to design continuum controllers may lead to problematic ``spillover effects" when the solutions to discretized formulations do not appropriately converge to the true continuum solution \cite{morris2020controller, gibson1991approximation} leading to excitement of problematic modes or even instability. This motivates the need for a development of an SLS framework in such continuous spatial domain settings.  

The work of \cite{conger2025convex} presents a significant step in this direction. Notably,  these results were restricted to finite-time horizons and discrete-time settings to enable tractable numerical solutions \cite{conger2025convex}, but 
\textit{it remains unclear whether an analogous implementation of the controller can be obtained without such discrete-time and finite-time-horizon assumptions, while preserving internal stability.} 
This work builds on the foundation in \cite{conger2025convex} to address this gap via an informative case study for control of the diffusion equation.  
The main results presented in this paper are as follows: 
\begin{itemize}
    \item The infinite time horizon LQR problem for the diffusion equation over the torus is equivalently formulated as an optimization over the closed-loop mappings, analogous to the finite-dimensional SLS framework (Theorem~\ref{thm:main}).
    \item This closed-loop formulation is solved analytically to recover the known LQR (Sections~\ref{sec:analytical}, \ref{sec:LQR}).
    \item An implementation of the optimal controller is constructed from the optimal closed-loop mappings, and is rigorously shown to preserve internal stability of the full infinite-dimensional system (Lemma~\ref{lem:Controller_Implementation}).
    \item A class of spatio-temporal constraints on the closed-loop mappings are shown to preserve convexity of the closed-loop design problem (Section~\ref{sec:constraints}). 
\end{itemize}

%We consider the continuous-time, infinite-time-horizon Linear Quadratic Regulator (LQR) problem for the diffusion equation over the torus \cite{troltzsch2010optimal, gibson1979riccati}, \cite[Ch. 7]{morris2020controller}. For this specific system, we demonstrate that an analytic expression for the well-known LQR gain can be obtained through a closed-loop design procedure analogous to the SLS framework in finite dimensions. Moreover, we \emph{derive an internally stable implementation of the controller as a feedback interconnection of the closed-loop mappings} derived from this design procedure. 
This work provides a step toward developing a general SLS-like framework for optimal controller design for distributed parameter systems in continuous-time and infinite time-horizon settings, which will complement the computational results and theory developed in \cite{conger2025convex}. Although the diffusion equation itself is not particularly interesting, it presents a proof of concept that can inform the use of similar methodology for other first-order in time linear, constant coefficient PDE systems.

%This paper is structured as follows. Section~\ref{sec:preliminaries} provides notation and mathematical preliminaries. Section~\ref{sec:setup} introduces the diffusion dynamics and the control design problem. Section~\ref{sec:CL} presents a novel formulation of the control design problem in terms of the systems ``closed-loop maps". This formulation is solved analytically in Section~\ref{sec:analytical}, and an \textit{ internally stable implementation} of the resulting controller is constructed from the optimal ``closed-loop maps". A comparison to the Riccati equation solution of the LQR problem is presented in Section~\ref{sec:LQR}. Finally, Section~\ref{sec:constraints} illustrates that spatio-temporal constraints can be imposed on the closed-loop maps while preserving convexity of the underlying controller design problem, thus enabling spatio-temporally constrained design of controller implementations. 

\section{Notation \& Mathematical Preliminaries} \label{sec:preliminaries}

Let $\mathbb{T}$ denote the unit circle, $\mathbb{R}$ the real numbers, $\mathbb{C}$ the complex numbers, and $\mathbb{Z}$ the integers.
For $\mathcal{X}, \mathcal{U}$ two Hilbert spaces, $\mathcal{L}(\mathcal{X}, \mathcal{U})$ denotes the
set of linear operators from $\mathcal{X}$ to $\mathcal{U}$, and  $\mathcal{L}(\mathcal{X}) = \mathcal{L}(\mathcal{X}, \mathcal{X}).$ An operator $B \in \mathcal{L}(\mathcal{X}, \mathcal{U})$ is \emph{bounded} if $\sup_{u \in \mathcal{U}} \frac{\|Bu\|_{\mathcal{X}}}{\|u\|_{\mathcal{U}}} < \infty$. 
$C^{\dagger} \in \mathcal{L}(\mathcal{U}, \mathcal{X})$ denotes the adjoint of $C \in \mathcal{L}(\mathcal{X}, \mathcal{U})$.
 Two Hilbert spaces of interest are the space $L^2(\mathbb{T})$ equipped with inner product  
    $\left< f, g \right>_{L^2(\mathbb{T})} = \int_{\theta \in \mathbb{T}} f^*(\theta) g(\theta) d \theta,$
and the space
$\ell^2(\mathbb{Z})$ of sequences equipped with inner product 
 $
    \big<\hat{f}, \hat{g} \big>_{\ell^2(\mathbb{Z})} = \sum_{\kappa \in \mathbb{Z}} \hat{f}^{*}(\kappa)\hat{g}(\kappa).
$

We consider \emph{spatio-temporal signals} over the domain $[0, \infty) \times \mathbb{T}, $ which we denote by lower-case letters e.g., 
\be 
    \psi(t, \theta), ~~ t \in [0, \infty), ~ \theta \in \mathbb{T}.
\ee 
We assume that such spatio-temporal signals $\psi(t, \cdot) \in L^2(\mathbb{T})$ for each $t$ and use  bold-face lower case letters to denote the corresponding $L^2(\mathbb{T})$ valued signals, e.g., 
$\bpsi(t) := \psi(t, \cdot).$

The Fourier transform, $\mathcal{F}$, of elements of $L^2(\mathbb{T}):$
\be 
    \hat{f}(j \kappa):= (\mathcal{F}f)(j \kappa): = \frac{1}{2 \pi} \int_{\theta \in \mathbb{T}} f(\theta) e^{-j \kappa \theta}d \theta, 
\ee 
 can be applied to the spatial component of the spatio-temporal signal $\psi(\cdot, \cdot)$: 
\be 
    \hat{\psi}(t, j \kappa) = \frac{1}{2 \pi} \int_{\theta \in \mathbb{T}} \psi(t, \theta) e^{-j \kappa \theta}d \theta,
\ee 
such that for each $t,$ $\boldsymbol{\psi}(t, \cdot) \in \ell^2(\mathbb{Z}).$
%Recall that $\mathcal{F}$ is an isometric isomorphism to  $\ell^2(\mathbb{Z})$.

Letting $+$ denote the addition operator on the group $\mathbb{T}, $ for each $\xi \in \mathbb{T}$, define the translation operator $T_{\xi} \in \mathcal{L}(L^2(\mathbb{T})):$
    $(T_{\xi}f)(\theta) := f(\theta +  \xi).$
An operator $B \in \mathcal{L}(L^2(\mathbb{T}))$ is \emph{translation invariant} if it commutes with every translation operator, i.e. $T_{\xi} B = B T_{\xi}$ for all $\xi \in \mathbb{T}.$ Such an operator can be represented in the Fourier domain as a \emph{multiplication operator}, $B = \mathcal{F}^{-1} M_{\hat{B}} \mathcal{F}$ with $M_{\hat{B}} \in \mathcal{L}(\ell^2(\mathbb{Z}))$ such that
\be 
    (M_{\hat{B}} \hat{f})( \kappa ) = \hat{B}( \kappa) \hat{f}( \kappa).
\ee 
We refer to the measurable function $\hat{B}_{\kappa} = \hat{B}(\kappa)$ as the \emph{symbol} of the operator $M_{\hat{B}}$.

\section{Problem Set Up} \label{sec:setup}
The controller and plant of interest are linear, time- and spatially-invariant distributed parameter systems over the spatial domain $\mathbb{T}$. These systems can be specified in abstract ODE form \cite{curtain2020introduction, morris2020controller} as  
\be \begin{aligned} \label{eq:generalODE}
    \tfrac{d}{dt} {\bpsi}(t) &= A \bpsi(t) + B (\bu + \bw)(t)\\
     \bz(t) &= C \bpsi(t) + D \bu (t),
\end{aligned} \ee 
where $A, B, C, D$ are spatially-invariant operators, $A$ generates a $C_0-$semigroup $\{T_A(t) \}$ on $\mathcal{D}(A) \subset L^2(\mathbb{T})$, $B \in \mathcal{L}(L^2(\mathbb{T}))$, and $C, D \in \mathcal{L}(L^2(\mathbb{T}), L^2(\mathbb{T}) \oplus L^2(\mathbb{T}))$. $\boldsymbol{u,w, \psi}$ and $\bz$ denote the control signal, external disturbance, state, and performance output signal, respectively.  We denote the system \eqref{eq:generalODE} by $\Sigma(A,B,C, D)$.
A linear, spatially-invariant state feedback controller for this system takes the form 
\be 
    {\boldsymbol{u}}(t) = K \boldsymbol{\psi}(t),
\ee 
where $K \in \mathcal{L}(L^2(\mathbb{T}))$ is a translation invariant operator. We denote the system resulting from interconnection of \eqref{eq:generalODE} with controller $K$ as $\mathbb{F}(P;K)$; this system is also spatially invariant and is described by 
\be \begin{aligned}
    \tfrac{d}{dt} {\bpsi}(t) &= (A+BK) \bpsi(t) + B {\bf w}(t)\\
     {\bf z}(t) &= (C+DK) \bpsi(t).
\end{aligned} \ee 
The $C_0$-semigroup $\{T_{A_{\rm c.l.}}(t)\}$ generated by $A_{\rm c.l.} := A+BK$ is 
\emph{exponentially stable} if for some constants $m, \eta > 0$,
\be 
    \| T_{A_{\rm c.l.}}(t)\|_{L^2(\mathbb{T}) \rightarrow L^2(\mathbb{T})} \le m e^{-\eta t}, ~~ \text{for all}~ t \ge 0.
\ee  
%{The system is \emph{output stable} if $CT_A(\cdot)$ is a bounded map from $L^2(\mathbb{T})$ to $L^2([0, \infty), L^2(\mathbb{T}))$.}

\begin{definition} The \emph{infinite horizon observability Gramian}, $L_O$, of the pair $(A_{\rm c.l.},C)$ is \cite{curtain2020introduction}
\be \label{eq:LO}
    L_O %= \mathcal{C}_{\infty}^{\dagger} \mathcal{C}_{\infty} 
    = \int_0^{\infty}(T_{A_{\rm c.l.}}(t))^{\dagger} C^{\dagger} C T_{A_{\rm c.l.}}(t) dt,
\ee 
when this integral converges absolutely in operator norm. A sufficient condition for $L_O$ to be well-defined is that $C$ is bounded and that $\{T_{A_{\rm c.l.}}(t)\}$ is exponentially stable. When $B$ is bounded and $L_O$ is trace class, the $\mathcal{H}_2$ \emph{norm}\footnote{Since we assumed $\bw(t)\in L^2(\mathbb{T})$, we can equivalently interpret $\|\cdot\|_{\mathcal{H}_2}$ as the $L^2$ norm of the impulse response as outlined in \cite[Ch. 5.1]{morris2020controller}.} of the system $\mathbb{F}(P;K) = \Sigma(A_{\rm c.l.}, B, (C+DK), 0)$ is defined by 
    \be \label{eq:H2_norm}
        \| \Sigma(A_{\rm c.l.},B,(C+DK),0) \|_{\mathcal{H}_2} = \tfrac{1}{2 \pi} {\rm Tr}(B^\dagger L_O B).
    \ee 
When $ A_{\rm c.l.} $ generates an exponentially stable $C_0$-semigroup 
and $\| \Sigma(A_{\rm c.l.},B,(C+DK),0) \|_{\mathcal{H}_2}<\infty$, 
we say that the controller $K$ is \emph{internally stabilizing}.
\end{definition}

%\be \begin{aligned}
%   & \mathcal{F}^{-1} \left( C(sI-A)^{-1} B + D\right) \mathcal{F} = M_{\hat{G}(s)}\\
   %& \hat{G}_{\kappa}(s) =\hat{C}_{\kappa} (sI - \hat{A}_{\kappa})^{-1}\hat{B}_{\kappa} + \hat{D}_{\kappa}.
%\end{aligned} \ee 
%and we say that the transfer function is spatially-invariant.
\subsection{Diffusion Dynamics}
We focus on a specific plant  - the diffusion equation over $\mathbb{T}$ with fully distributed actuation $u$ and disturbance $w$, described by the partial differential equation (PDE) 
\be \label{eq:olPDE}
    \partial_t \psi(t, \theta) = \alpha \partial_{\theta}^2\psi(t, \theta) + u(t, \theta) + w(t, \theta), ~t \in [0, \infty), ~ \theta \in \mathbb{T},
\ee
where $u(t,\cdot),w(t,\cdot),\psi(t,\cdot)\in L^2(\mathbb{T})$ and $\alpha > 0$ is the constant of diffusivity.  The 
 spatio-temporal signal
 \be \label{eq:output_signal}
    {z}(t, \theta) = \lba{c} \gamma {\psi}(t, \theta) \\ { u}(t, \theta) \ear, ~z(t, \cdot) \in L^2(\mathbb{T}) \oplus L^2(\mathbb{T})
\ee
is a performance output that provides a uniform penalty across the spatial domain
with $\gamma >0$ determining a tradeoff between cost of state and cost of control.
The interest in this system is not for its particular physical significance, but rather to serve as a case study for the proposed framework. Restricting to the torus is equivalent to enforcing periodic boundary conditions and will admit analytic solutions via Fourier analysis approaches \cite{bamieh2002distributed}.

The dynamics \eqref{eq:olPDE}-\eqref{eq:output_signal} can be written in the form \eqref{eq:generalODE} as 
\be \begin{aligned} \label{eq:olODE}
    \tfrac{d}{dt} {\bpsi}(t) &= A \bpsi(t) +  {\bu}(t) + {\bw}(t)\\
     {\bz}(t) &= C \bpsi(t) + D {\bu}(t),
\end{aligned} \ee 
where $A \in \mathcal{L}(L^2(\mathbb{T}))$ is defined by 
$(A \psi)(t, \theta) = \alpha \partial_{\theta}^2\psi(t,\theta)$
and generates a $C_0-$ semigroup $\{T_A(t)\}$ on 
\be \label{eq:domain}
   \{ f \in L^2(\mathbb{T});~ f', f'' \in L^2(\mathbb{\mathbb{T}}) \} =: \mathcal{D}(A) \subset L^2(\mathbb{T}),
\ee 
with derivatives, $f', f''$, in \eqref{eq:domain} interpreted weakly, 
and $C,D \in \mathcal{L}(L^2(\mathbb{T}),L^2(\mathbb{T}) \oplus L^2(\mathbb{T}))$ are  the bounded operators 
$$ \begin{aligned}%\label{eq:CD}
   (C\psi)(t, \theta) = \lba{c} \gamma \psi(t, \theta) \\ 0 \ear,~ (Du)(t, \theta) = \lba{c} 0 \\ u(t, \theta)\ear .
\end{aligned} $$
$\Sigma(A,I,C,D)$ is a spatially-invariant system, equivalently characterized by its transformation under $\mathcal{F}$ as 
\be \begin{aligned} \label{eq:open_loop_kappa}
    &\tfrac{d}{dt} \hat{\psi}(t,  \kappa) = %\underbrace{
    - \alpha \kappa^2
    %}_{\hat{A}_{\kappa}} 
    \hat{\psi}(t,  \kappa) + \hat{u}(t,  \kappa) + \hat{w}(t,  \kappa), \\
    &\hat { z} (t,  \kappa)  = \underbrace{\lba{c} \gamma \\ 0 \ear }_{\hat{C}_{\kappa}}\hat{\psi}(t, \kappa) + \underbrace{\lba{c} 0 \\ 1 \ear}_{\hat{D}_{\kappa}} \hat{u}(t, \kappa), ~ t \in [0, \infty), ~\kappa \in \mathbb{Z}.
\end{aligned} \ee 
Note that for each fixed $\kappa \in \mathbb{Z},$ \eqref{eq:open_loop_kappa} is a system with finitely many inputs and outputs.

\subsection{Optimal Control Formulations}
We consider the Linear Quadratic Regulator (LQR) problem for the diffusion equation \eqref{eq:olODE} \cite{banks1984linear, curtain2020introduction}: 
\be \begin{aligned} \label{eq:LQR}
    &\inf_{u}~ \int_{0}^{\infty} \left< \bpsi(t), \gamma^2 \bpsi(t) \right>_{L^2(\mathbb{T})} + \left<{\bu}(t), {\bu}(t) \right>_{L^2(\mathbb{T})} dt\\
    &~{\rm s.t.}~~ \partial_t \psi(t, \theta) = \alpha \partial_{\theta}^2\psi(t, \theta) + u(t, \theta), ~ {\bpsi}(0) = \psi_0.
    %&~~~~~~~{\bf u}(t) = F \bpsi(t) ~\text{internally stabilizing,}
\end{aligned}\ee 
 The solution $K^{\rm opt}$ to \eqref{eq:LQR} is independent of initial condition $\psi_0 \in L^2(\mathbb{T})$ and it will be shown in Section~\ref{sec:LQR} that $K^{\rm opt}$ is the solution to the following  $\mathcal{H}_2$ problem:  

\be \begin{aligned} \label{eq:H2_formulation}
   K^{\rm opt} =  &\argmin_K ~~ \| \mathbb{F}(P;K) \|_{\mathcal{H}_2}\\
    &~~~~{\rm s.t.}~~~~ K ~\text{spatially-invariant}\\
    &~~~~~~~~~~~~K ~\text{internally stabilizing}.
\end{aligned}\ee 
\begin{remark}
   In some cases \eqref{eq:LQR} and \eqref{eq:H2_formulation} result in the same $K^{\rm opt}$, but this does not always hold for PDE systems. E.g., if operation by $K$ requires convolution with a delta distribution, the impulse response of the closed-loop system will not be $L^2(\mathbb{T})$; thus, the $\mathcal{H}_2$ problem is not well-posed, although the LQR problem is well-posed with finite cost $\forall ~\psi_0\in \mathcal{D}(A)$. 
\end{remark}

%\begin{definition} We say that the control policy \eqref{eq:K_TF} \emph{internally stabilizes} the system \eqref{eq:olODE} if the transfer function $\mathbb{F}(P;K)(s)$ is exponentially stable
%$C_0$-semigroup $\{T_{A_{\rm c.l.}}(t) \}$ and $\|\mathbb{F}(P;K)\|_{\mathcal{H}_2}< \infty$.
% \end{definition}

%As in the finite-dimensional setting, this $\mathcal{H}_2$ problem will recover the solution to the Linear Quadratic Regulator (LQR) problem for the diffusion equation 

\section{Closed-Loop Formulation of $\mathcal{H}_2$ Problem} \label{sec:CL}
In this section, we demonstrate that the controller design problem \eqref{eq:H2_formulation} can be formulated as an optimization problem over certain closed-loop mappings and admits analytical solutions in the spatial frequency domain. We will illustrate that this formulation allows for spatio-temporal constraints to be imposed in a convex manner in Section~\ref{sec:constraints}. % , thus providing a first step toward a tractable formulation of optimal controller design problems for distributed parameter systems with spatio-temporal constraints in the continuous-time, infinite time horizon setting.

%Our main result provides an optimization problem equivalent to \eqref{eq:H2_formulation} that admits a class of spatio-temporal structural constraints while preserving convexity. 
%We begin by defining certain ``closed-loop maps" - operator-valued transfer functions that are analogous to the finite-dimensional transfer functions optimized over in the System Level Synthesis methodology \cite{wang2019system}. 
\subsection{Transfer Functions} 
To introduce closed-loop mappings, we use operator-valued transfer function representations of systems \cite[Ch.7]{curtain2020introduction} to denote mappings between temporally Laplace transformed variables.  Under certain technical conditions \cite[Ch. 7]{curtain2020introduction}, for some positive constant $\beta$, the Laplace transform of a signal ${\bpsi}$, ${\bPsi}: \{s \in \mathbb{C}; ~ {\rm Re}(s) \ge \beta \} \rightarrow L^2(\mathbb{T})$, is well-defined by  
   $ \bPsi(s) := \int_0^{\infty} e^{-st} \bpsi(t) dt.$ First let $\Sigma(A,B,C,D)$ be a general system \eqref{eq:generalODE}. Let $\rho(A)$ denote the resolvent set and $\mathcal{D}(A)\subset L^2(\mathbb{T})$ denote the domain of the operator $A$. For $s \in \rho(A)$, the \emph{resolvent operator} $ (sI -A)^{-1}$ is well-defined on $\mathcal{D}(A) $ by
\begin{equation} 
  %& R(s;A) = (sI-A)^{-1}: D(A) \rightarrow L^2(\mathbb{T})\\
   %& \psi_0 \in\mathcal{D}(A) \mapsto (sI-A)^{-1} \psi_0:=   \int_0^{\infty} e^{-st}T_A(t) \psi_0 dt.
   (sI-A)^{-1}\psi_0:=   \int_0^{\infty} e^{-st}T_A(t) \psi_0 dt,~\psi_0\in\mathcal{D}(A).
\end{equation}
$\rho(A)$ will contain a set of the form $\{ s\ \in \mathbb{C}; ~ {\rm Re}(s) \ge \omega\}$, where $\omega$ defines a growth bound on the semigroup $\{T_A(t)\}$. When this semigroup is exponentially stable, one can take $\omega = 0 $ \cite[Ch. 2]{curtain2020introduction}.
Then the operator-valued \emph{transfer function} representation of a system \eqref{eq:generalODE} $\Sigma(A,B,C,D)$ is given by 
\be 
     C(sI-A)^{-1} B + D: \mathcal{D}(A) \rightarrow L^2(\mathbb{T}),~~ s \in \rho(A).
\ee 
 If $A$ generates an exponentially stable $C_0-$semigroup, then we say that the transfer function representation of $\Sigma(A,B,C,D)$ is \emph{exponentially stable}. If $D\equiv0$, the transfer function representation of $\Sigma(A,B,C,D)$ is \emph{strictly proper}. 
 
Now let $\Sigma(A,B,C,D)$ describe the diffusion equation \eqref{eq:olODE}. We denote the operator-valued transfer function from inputs $\begin{bmatrix} \boldsymbol{w} & \boldsymbol{u} \end{bmatrix}^\top $ to $\begin{bmatrix}
    \boldsymbol{z} & \boldsymbol{\psi}
\end{bmatrix}^\top$ for \eqref{eq:olODE} by  $P(s)$,
which is well-defined on { $\{s \in \mathbb{C}; {\rm Re}(s) > 0 \} \subset \rho(A)$} as $0$ defines a growth bound on the semigroup $\{T_A(t) \}.$ 
Since \eqref{eq:olODE} is a spatially invariant system, $P(s)$ can be expressed as a multiplication operator, $M_{\hat{P}(s)},$ over spatial frequency \cite[Ch.7]{curtain2020introduction}, and \eqref{eq:open_loop_kappa} can be represented by the transfer function
\be \begin{aligned}\label{eq:transfer_function_state}
   % &,\\& 
   \hat{P}_{\kappa}(s)  = \begin{bmatrix}
   \hat{C}_{\kappa} (s+ \alpha \kappa)^{-1} & \hat{C}_k (s+ \alpha \kappa)^{-1} + \hat{D}_{\kappa} \\
    (s+ \alpha \kappa)^{-1} &  (s+ \alpha \kappa)^{-1} 
   \end{bmatrix},
\end{aligned} \ee 
which maps $\begin{bmatrix} \hat{w}_{\kappa}(s) & \hat{u}_{\kappa}(s) \end{bmatrix}^\top$ to $\begin{bmatrix} \hat{z}_{\kappa}(s) & \hat{\psi}_{\kappa}(s) \end{bmatrix}^\top$. 

We consider a (possibly dynamic) state feedback control policy defined by a spatially invariant transfer function
\be \label{eq:K_TF}
    {\bf U}(s) = K(s) \bPsi(s),
\ee
defined on some set $\{s \in \mathbb{C};~ { \rm Re}(s) \ge  \omega\}$.
The feedback interconnection of $P(s)$ and \eqref{eq:K_TF}, denoted $\mathbb{F}(P; K)(s)$, is a well-defined transfer function on $\{s \in \mathbb{C}; {\rm Re}(s) > 0 \} \cap \{s \in \mathbb{C};~ { \rm Re}(s) \ge  \omega\}$ %$\{\color{red}s \in \mathbb{C}; ~ {\rm Re}(s) > \max\{0, \omega\} \}$ 
\cite[Ex. 7.3]{curtain2020introduction} that maps $\boldsymbol{w}(s)$ to $\boldsymbol{z}(s)$. 
 
 \subsection{Closed-Loop Mappings}
We isolate the closed-loop transfer functions from $\bw$ to $\bpsi$ and $\bu$ in $\mathbb{F}(P; K)(s)$ to obtain: 
\be \begin{aligned} \label{eq:Phi}
    \bPsi(s) &= (I - (sI-A)^{-1}K(s))^{-1} (sI-A)^{-1}{\bf W}(s) \\
    & = (sI - A - K(s))^{-1} {\bf W}(s)\\ 
    & =: \Phi^{\psi}(s) {\bf W}(s)\\
     {\bf U}(s) &= K(s) (sI - A - K(s))^{-1} {\bf W}(s)\\
     & =: \Phi^{u}(s) {\bf W}(s).
\end{aligned} \ee 
The operator-valued transfer functions $\Phi^{\psi}(s), \Phi^u(s) \in \mathcal{L}(L^2(\mathbb{T}))$ for each $s\in \mathbb{C}$ at which they are defined and are analogous to the ``system-level maps" introduced for finite dimensional systems in \cite{wang2019system}. 
The closed-loop transfer function from $\bw$ to $\bz$ can be written in terms of $\Phi^{\psi}(s), \Phi^u(s)$ as 
\be\begin{aligned} \label{eq:Phi_to_F}
     {\bf Z}(s) & = \mathbb{F}(P;K)(s) {\bf W}(s)\\
     &=   \left( C \Phi^{\psi}(s) + D \Phi^u(s) \right){\bf W}(s).
\end{aligned}\ee 
We parameterize the set of internally stabilizing, time- and spatially-invariant state feedback controllers for \eqref{eq:olODE} in terms of 
 $\Phi^{\psi}(s)$ and $ \Phi^u(s)$.

\begin{lemma} \label{lem:KStabilizing_PhiStable}
    Let $K$ be a spatially-invariant and internally stabilizing controller for \eqref{eq:olODE}. Then the operator-valued transfer functions %$\Phi^{\psi}(s)$ and $\Phi^u(s)$ 
    \eqref{eq:Phi} are spatially-invariant, strictly proper, well-defined on $\{s \in \mathbb{C};~ {\rm Re}(s) \ge 0 \},$ and satisfy 
    \be \begin{aligned} \label{eq:Hinf_condition}
        &\sup_{{\rm Re}(s) \ge 0} ~ \| \Phi^{\psi}(s) \|_{L^2(\mathbb{T}) \rightarrow L^2(\mathbb{T})} < \infty, \\
       & \sup_{{\rm Re}(s) \ge 0} ~ \| \Phi^{u}(s) \|_{L^2(\mathbb{T}) \rightarrow L^2(\mathbb{T})} < \infty.
    \end{aligned} \ee 
    Condition \eqref{eq:Hinf_condition} can be equivalently verified pointwise in the spatial frequency domain as: 
    \be \begin{aligned} \label{eq:Hinf_frequency}
        &\sup_{{\rm Re}(s) \ge 0} \sup_{\kappa \in \mathbb{Z}} ~\left| \hat{\Phi}^{\psi}_{\kappa}(s)  \right| < \infty , ~ \sup_{{\rm Re}(s) \ge 0} \sup_{\kappa \in \mathbb{Z}} ~\left| \hat{\Phi}^{u}_{\kappa}(s)  \right| < \infty.
    \end{aligned} \ee 
    \end{lemma}

\begin{proof} See Appendix~\ref{app:Phi_in_Hinf}.
\end{proof}

Given the analogy to the definition for finite dimensional systems, we say that an operator-valued transfer function is an element of $\mathcal{H}_{\infty}$ if it satisfies the condition \eqref{eq:Hinf_condition}. %By the Plancherel Theorem, 
%$\sup_{{\rm Re}(s) \ge 0} ~ \| \Phi^{\psi}(s) \|_{L^2(\mathbb{T}) \rightarrow L^2(\mathbb{T})}  = \sup_{{\rm Re}(s) \ge 0} ~ \| M_{\hat{\Phi}^{\psi}(s)} \|_{\ell^2(\mathbb{Z}) \rightarrow \ell^2(\mathbb{Z})},$ and similarly for $\Phi^u.$ x

\begin{lemma}\label{lem:op_constraint}
    Let %the operator-valued transfer functions 
    $\Phi^{\psi}(s)$, $\Phi^u(s) \in \mathcal{L}(L^2(\mathbb{T}))$ be defined by \eqref{eq:Phi} for some spatially-invariant and internally stabilizing controller $K(s).$ Then, for all $\eta, \vartheta \in L^2(\mathbb{T})$ and $s \in \{s \in \mathbb{C}; {\rm Re}(s) \ge 0 \}$,
    \be \begin{aligned} \label{eq:affine_innerProduct}
       & \left<\eta, \left((s+1)I- (A+I)\right)\Phi^{\psi}(s) \vartheta \right>_{L^2(\mathbb{T})}\\&\hspace{1.
       in}- \left< \eta, \Phi^u(s) \vartheta \right>_{L^2(\mathbb{T})} 
        = \left< \eta, \vartheta \right>_ {L^2(\mathbb{T})}.
    \end{aligned} \ee 
    Equivalently,
     \be \label{eq:affine_pointwise}
    \left( (s+1) - (- \alpha \kappa^2 + 1) \right) \hat{\Phi}^{\psi}_{\kappa}(s) - \hat{\Phi}_{\kappa}^u(s) = 1, ~ \forall \kappa \in \mathbb{Z}.
 \ee 
\end{lemma}
\begin{proof}
    See Appendix~\ref{app:op_constraint}.
\end{proof}
We can now state our main result: 
\begin{theorem}\label{thm:main}
    The solution to the $\mathcal{H}_2$ design problem for the diffusion equation \eqref{eq:H2_formulation} can be recovered from the solution to:
\be \begin{aligned} \label{eq:CL_formulation}
    &\argmin_{\Phi^{\psi}, \Phi^u \in \mathcal{H}_{\infty}}~ \left\| C \Phi^{\psi}  + D \Phi^u\right\|_{\mathcal{H}_2} \\
    &~~~~~{\rm s.t.}~~~~~~\Phi^{\psi}, \Phi^u \text{ spatially invariant, strictly proper}\\
    & ~~~~~~~~~~~~~~~ \Phi^{\psi}, \Phi^u \text{ satisfy } \eqref{eq:affine_innerProduct}.
\end{aligned}
\ee 
\end{theorem}

Lemmas \ref{lem:KStabilizing_PhiStable}-\ref{lem:op_constraint} show that $K$ spatially invariant and internally stabilizing implies that $\Phi^{\psi}, \Phi^u \in \mathcal{H}_{\infty}$ and are spatially invariant, strictly proper, and satisfy the affine constraint \eqref{eq:affine_innerProduct}. To complete the proof of Theorem \ref{thm:main} it remains to show that any such $\Phi^{\psi}, \Phi^u$ can be implemented as an internally stabilizing controller. We will prove this in Lemma \ref{lem:Controller_Implementation} and discuss controller implementation in Section \ref{sec:imp}.

\subsection{Analytical Solution to Closed-Loop Formulation of $\mathcal{H}_2$ Problem} \label{sec:analytical}

Spatial invariance  allow us to solve \eqref{eq:CL_formulation} in the spatial frequency domain. By the Plancherel theorem, formulation \eqref{eq:CL_formulation} of the $\mathcal{H}_2$ design problem for the diffusion equation is equivalent to the optimization problem: 
    \be \begin{aligned} \label{eq:CLOpt_Frequency}
            &\argmin_{{\hat{\Phi}^{\psi}(s)}, {\hat{\Phi}^u(s)}} ~ \left\| M_{\hat{C}} M_{\hat{\Phi}^{\psi}(s)} + M_{\hat{D}} M_{\hat{\Phi}^u(s)}\right\|_{\mathcal{H}_2}\\
        & ~~~~~{\rm s.t.}~~~~~~\hat{\Phi}^\psi,\hat{\Phi}^u \text{ strictly~proper,~spatially~invariant}, \\
        & ~~~~~~~~~~~~~~~\hat{\Phi}^\psi,\hat{\Phi}^u \text{ satisfy \eqref{eq:Hinf_frequency} and \eqref{eq:affine_pointwise}.}
%        & ~~~~~~~ \sup_{{\rm Re}(s) \ge 0} \sup_{\kappa \in \mathbb{Z}} ~\left| {{\hat{\Phi}^{\psi}_{\kappa}(s) }} \right| < \infty,\\
       % & ~~~~~~~\sup_{{\rm Re}(s) \ge 0} \sup_{\kappa \in \mathbb{Z}} ~\left| \hat{\Phi}^{u}_{\kappa}(s)  \right| < \infty
    \end{aligned}
    \ee 
Note that the cost and constraints of \eqref{eq:CLOpt_Frequency} ``decouple" in the spatial frequency variable, so the solution to \eqref{eq:CLOpt_Frequency} can be obtained by solving a parameterized family of $\mathcal{H}_2$ problems for finite-dimensional systems parameterized by $\kappa \in \mathbb{Z}$. Let ${\mathcal {RH}}_2$ denote the set of rational transfer functions with strictly proper entries and no poles in the closed right half plane \cite[Ch. 3]{dullerud2013course}. Then the symbol of the minimizer of \eqref{eq:CLOpt_Frequency} at frequency $\kappa$ is the solution to the family of problems:
\begin{subequations}\label{eq:CLOpt_kappa} \begin{align} 
    &\argmin_{\hat{\Phi}^{\psi}_{\kappa}, \hat{\Phi}^u_{\kappa}\in \mathcal{RH}_{2}} \left\| \lba{cc} \gamma & 0  \\ 0 & 1   \ear \lba{c} \hat{\Phi}^{\psi}_{\kappa} \\ \hat{\Phi}^u_{\kappa}\ear \right\|_{\mathcal{H}_2}\\
    &~~{\rm s.t.} \left( (s+1) - (- \alpha \kappa^2 + 1) \right) \hat{\Phi}^{\psi}_{\kappa}(s) - \hat{\Phi}_{\kappa}^u(s) = 1, \label{eq:implicit_constraint}
\end{align} 
\end{subequations} 
provided that the family of solutions $\{\hat{\Phi}^{\psi}_{\kappa}\}, ~ \{\hat{\Phi}^{u}_{\kappa}\}$ satisfy the condition \eqref{eq:Hinf_frequency} which ensures $\Phi^{\psi}, \Phi^u \in \mathcal{H}_{\infty}$.  

For each fixed $\kappa \in \mathbb{Z},$ \eqref{eq:CLOpt_kappa} is a standard $\mathcal{H}_2$ design problem for a finite-dimensional system with an implicit affine subspace constraint \eqref{eq:implicit_constraint}, which we can solve analytically. 
\begin{proposition}\label{prop:projection}
    The optimal solution $\hat{\Phi}^{\psi}_{\kappa}, ~\hat{\Phi}^u_{\kappa}$ to the problem \eqref{eq:CLOpt_kappa}, with $\gamma$ defined in \eqref{eq:output_signal},  is given by 
    \be \begin{aligned} \label{eq:opt_Phi}
        \hat{\Phi}^{\psi}_{\kappa}(s) &= \left(s + \sqrt{\gamma^2 + \alpha^2 \kappa^4}\right)^{-1}, \\
        \hat{\Phi}^{u}_{\kappa}(s) &= \frac{\alpha \kappa^2
        - \sqrt{\gamma^2 + \alpha^2 \kappa^4}}{s + \sqrt{\gamma^2 + \alpha^2 \kappa^4}}.
    \end{aligned}\ee 
\end{proposition}

\begin{proof}
    See Appendix~\ref{app:projection}.
\end{proof}

\section{internally stable implementation of Controller via Closed-Loop Solution} \label{sec:imp}
\begin{figure}[h!]
\centering\includegraphics[width=.75\linewidth]{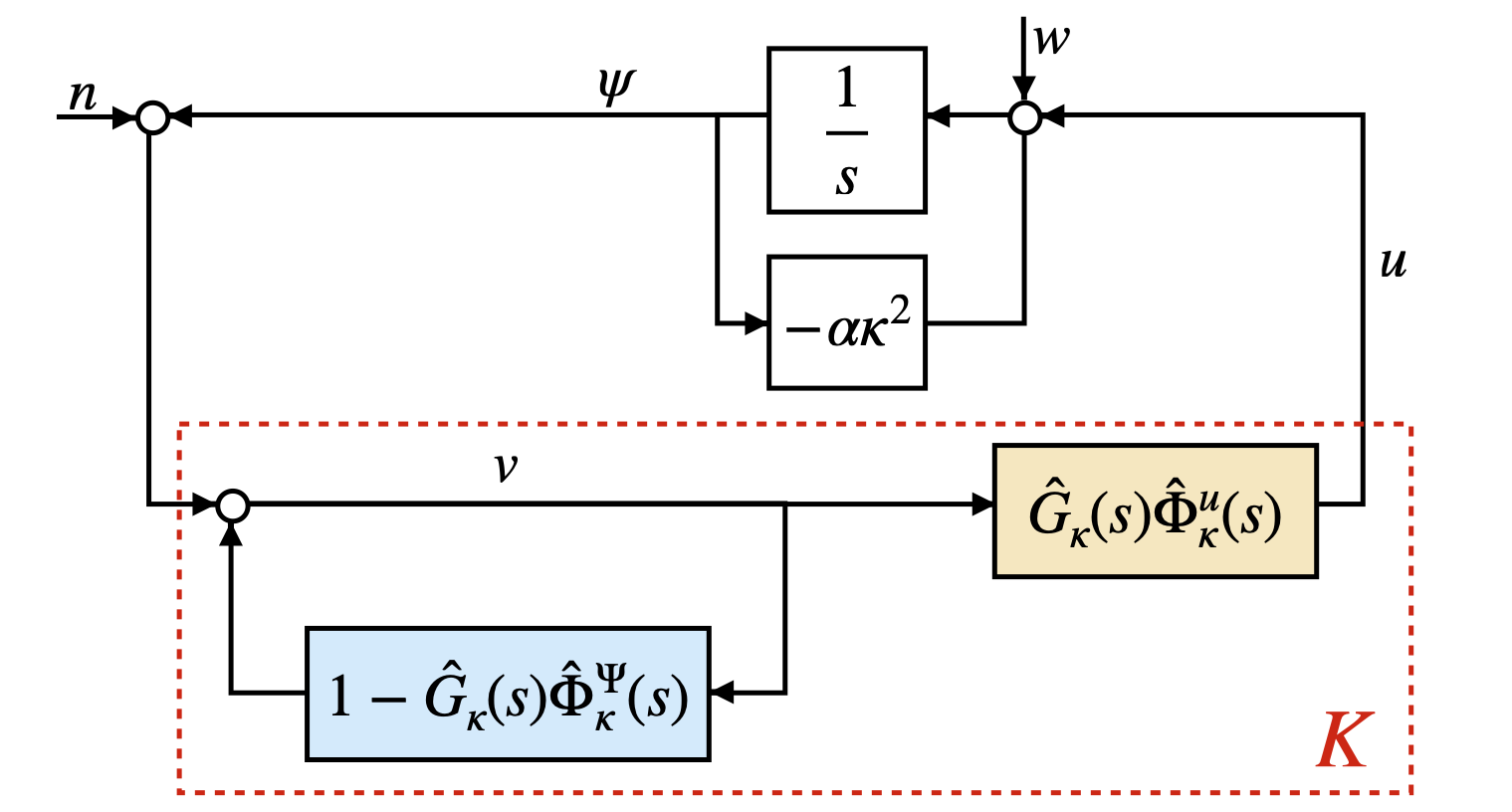}
    \caption{Block diagram of the dynamic controller implementation \eqref{eq:control_implementation} with noisy measurements and exogenous disturbance: $w$ is additive noise to the state $\psi$, and $w$ is additive noise to the control signal $u$.}
    \label{fig:controller_implementation_diagram}
\end{figure}

We now complete the proof of Theorem \ref{thm:main} by showing that $\Phi^{\psi}, \Phi^u$ satisfying the constraints of \eqref{eq:CL_formulation} can be used to construct an internally stable implementation of the controller, which thus preserves the structural properties of these maps. In order for an implementation of the controller to be well posed and preserve the structure of $\Phi^{\psi}$ and $\Phi^u$, we introduce an additional transfer function $G$ that ensures internal states of the controller remain bounded. 
\begin{lemma} \label{lem:Controller_Implementation}
    Let ${\Phi}^{\psi},{\Phi}^u \in \mathcal{H}_{\infty} \subset \mathcal{L}(L^2(\mathbb{T}))$ be spatially-invariant, strictly proper, and satisfy \eqref{eq:affine_pointwise}. Select %$\hat{g}_{\kappa}(s)$ such that {$\left(\hat{g}_{\kappa}(s)\right)^{-1} (sI-\hat{A}_{\kappa}) \in H_{\infty}$ }
    a spatially-invariant transfer function $G$ for which $G^{-1}(s)\in \mathcal{H}_\infty$, $\lim_{s \to \infty} s\hat G_\kappa^{-1}(s) = 1 ~\forall \kappa \in \mathbb{Z}$, and $ G^{-1}(s)(sI - A) \in \mathcal{H}_{\infty}$. Then,
    \be \begin{aligned} \label{eq:control_implementation}
        {\bf V}(s)& = {\bPsi}(s) + \left( I - G(s) {\Phi}^{\psi}(s)\right) {\bf V}(s)\\
        {\bf U}(s) & = G(s){\Phi}^u(s){\bf V}(s)
   \end{aligned} \ee
  %  \be \begin{aligned} \label{eq:control_implementation}
   %     \hat{v}(s, \kappa)& = \hat{\psi}(s, \kappa) + \left( 1 - \hat{g}_{\kappa}(s) \hat{\Phi}_{\kappa}^{\psi}(s)\right) \hat{v}(s, \kappa)\\
    %    \hat{u}(s,  \kappa) & = \hat{g}_{\kappa}(s)\hat{\Phi}^u_{\kappa}(s)\hat v(s, \kappa)
   %\end{aligned} \ee 
    defines an \emph{internally stable implementation} of the controller in the sense that the closed loop transfer matrices from the exogenous inputs $(\hat{n}_{\kappa}(s), \hat{w}_\kappa(s))$ to the internal signals $(\hat\psi_\kappa(s), \hat u_\kappa(s), \hat v_\kappa(s))$ belong to $\mathcal{H}_\infty$. This implementation is depicted in Figure~\ref{fig:controller_implementation_diagram}.
    The input-output map (transfer function) corresponding to \eqref{eq:control_implementation} is 
    \be \label{eq:recover_K}
        {\bf U} = \mathcal{F}^{-1}M_{\hat{K}(s)} \mathcal{F} {\bPsi},~~ \hat{K}_{\kappa}(s) = \tfrac{\hat{\Phi}^{u}_{\kappa}(s)}{\hat{\Phi}^{\psi}_{\kappa}(s)}.
    \ee
    %and the feedback interconnection of \eqref{eq:olODE} with the controller \eqref{eq:control_implementation} results in the closed-loop dynamics 
   % \be \label{eq:lemma_3_closed_loop_dynamics}
    %    {\bf Z}(s) = \left( C \Phi^{\psi}(s) + D \Phi^u(s) \right) {\bf W}(s).
   % \ee 
\end{lemma}

\begin{proof}
    See Appendix~\ref{app:controller_implementation}. 
\end{proof}

\begin{remark}
    In the finite-dimensional, discrete-time setting, the analogy of the system $G$ in \eqref{eq:control_implementation} is simply the one-time-step delay operator. More care must be taken to select this system in the infinite-dimensional setting due to the potential unboundedness of the operator $A$. In this case, the closed-loop map from additive noise to the state to the internal control signal $v$ is given by $(sI - A) G^{-1}(s)$ (see \eqref{eq:transfer_functions}), which may also be unbounded unless $G$ is carefully selected.
\end{remark}

\begin{figure}[h]
    \centering
    \includegraphics[width=\linewidth]{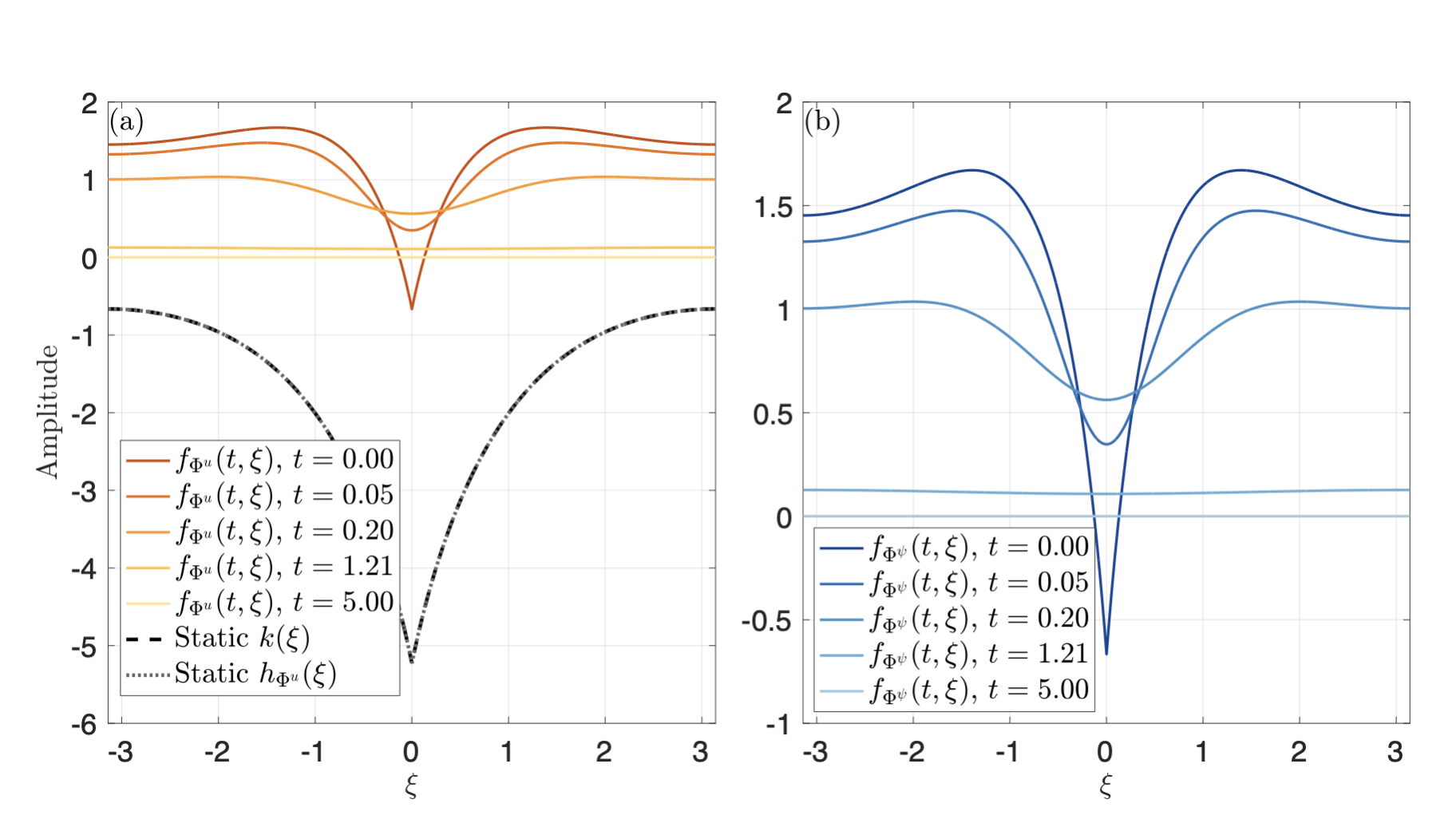}
    \caption{\small Spatio-temporal convolution kernels corresponding to the transfer functions in Fig. \ref{fig:controller_implementation_diagram}. Plot (a) shows the yellow block's static, $h_{\Phi^u}(\xi)$, and dynamic $f_{\Phi^u}(t, \xi)$ kernels. Plot (b) shows the blue block's dynamic kernel $f_{\Phi^\psi}(t, \xi)$.}
    \label{fig:kernel_comparisons}
\end{figure}

\label{sec:implementation}
We examine the structure of the controller implementation \eqref{eq:control_implementation} that results from the optimal solution \eqref{eq:opt_Phi}. We select $G(s) = \mathcal{F}^{-1}M_{\hat{G}(s)} \mathcal{F}$ where 
\be \begin{aligned} \label{eq:g}
    \hat{G}_{\kappa}(s) %&=  (1+\kappa^2) (s+1)\\
   %{\rm or}~ & = 
   = s + 1 + \kappa^2,
\end{aligned} \ee 
to satisfy all conditions of Lemma~\ref{lem:Controller_Implementation}. Note that this is one suitable choice of $G$, but it is not obvious how to ``best" select $G$ in general.
Then, \eqref{eq:control_implementation} is specified at frequency $\kappa$ as 
\be \begin{aligned} \label{eq:opt_implementation}
        \hat{v}(s, \kappa)& = \hat{\psi}(s, \kappa) + \left( 1 - \hat{G}_{\kappa}(s) \tfrac{1}{s + \sqrt{\gamma^2 + \alpha^2 \kappa^4}}\right) \hat{v}(s, \kappa)\\
        %& = \hat{\psi}(s, \kappa) + \left( 1 -  \tfrac{s + 1 + \kappa^2}{s + \sqrt{\gamma^2 + \alpha^2 \kappa^4}}\right) \hat{v}(s, \kappa)\\
        \hat{u}(s,  \kappa) & = \hat{G}_{\kappa}(s) \tfrac{\alpha \kappa^2 - \sqrt{\gamma^2 + \alpha^2 \kappa^4}}{s + \sqrt{\gamma^2 + \alpha^2 \kappa^4}}\hat{v}(s, \kappa),% ~= ~ \tfrac{(s+1+ \kappa)^2 \cdot \left( \alpha \kappa^2 - \sqrt{\gamma^2 + \alpha^2 \kappa^4} \right) }{s + \sqrt{\gamma^2 + \alpha^2 \kappa^4}}
\end{aligned} \ee 
where $\hat{G}$ is defined in \eqref{eq:g}. We note that the system $\mathcal{F}(I - G(s)\Phi^{\psi}(s))\mathcal{F}^{-1}$ with symbol $\left( 1 - \hat{G}_{\kappa}(s) \tfrac{1}{s + \sqrt{\gamma^2 + \alpha^2 \kappa^4}}\right)$ is strictly proper, so that the feedback loop in the control implementation in Fig. \ref{fig:controller_implementation_diagram} is well-defined.      Using the affine identity in \eqref{eq:affine_pointwise} together with routine calculations yields the following closed-loop transfer matrices from $(\hat n_\kappa(s), \hat w_\kappa(s))$ to $(\hat\psi_\kappa(s), \hat u_\kappa(s), \hat v_\kappa(s))$ in the controller implementation:
     \begin{align} \label{eq:transfer_functions}
         \begin{bmatrix}
             \hat\psi_\kappa\\
             \hat u_\kappa\\
             \hat v_\kappa
         \end{bmatrix} = \begin{bmatrix}
             \hat\Phi^u_\kappa(s) & \hat\Phi^\psi_\kappa(s)\\
             (s + \alpha\kappa^2)\hat\Phi^u_\kappa(s) & \hat\Phi^u_\kappa(s)\\
             (s + \alpha\kappa^2)\hat G^{-1}_\kappa(s)& \hat G^{-1}_\kappa(s)
         \end{bmatrix} \begin{bmatrix}
             \hat n_\kappa\\
             \hat w_\kappa\\
         \end{bmatrix},
     \end{align}
each of which is verified to be the symbol of a multiplication operator in $\mathcal{H}_{\infty}$ for our selected $\hat{G}.$ We numerically compute the spatio-temporal impulse responses (Green's functions) of the blocks in Fig. \ref{fig:controller_implementation_diagram} %using Matlab's \texttt{fft} function 
and plot their dynamic and static components over the spatial domain for various selections of time in Fig. \ref{fig:kernel_comparisons}. We note that the kernels of the dynamic components are ``more localized" for smaller time values and appear to spread out spatially as time increases.

\section{Equivalence to LQR Control} \label{sec:LQR}
We demonstrate that the solution to our closed-loop $\mathcal{H}_2$ design problem \eqref{eq:CL_formulation} recovers the solution to the standard infinite-horizon LQR problem for the diffusion equation.
\begin{proposition} \label{prop:LQR_equivalence} 
The standard LQR problem \eqref{eq:LQR} and the $\mathcal{H}_2$ problem of interest \eqref{eq:CL_formulation} have the same solution $K^{\rm opt}$.
\end{proposition}

\begin{proof} First, we compute the solution to the LQR formulation \eqref{eq:LQR} using standard approaches \cite{bamieh2002distributed, curtain2020introduction}.
We solve the operator-valued Riccati equation 
\be \begin{aligned} \label{eq:operator_riccati}
    \left< A\eta, F \vartheta \right>_{L^2(\mathbb{T})} +\left< F\eta, A \vartheta \right>_{L^2(\mathbb{T})} + \gamma^2 \left<  \eta,  \vartheta \right>_{L^2(\mathbb{T})} \\= \left< F \eta, F \vartheta \right>_{L^2(\mathbb{T})}.
\end{aligned} \ee 
The solution $F$ to \eqref{eq:operator_riccati} is known to be spatially-invariant \cite{bamieh2002distributed}, and the symbol of $\mathcal{F}F \mathcal{F}^{-1} = M_{\hat{F}}$ can be obtained by solving the family of scalar-valued algebraic Riccati equations (AREs) parameterized by $\kappa \in \mathbb{Z}:$
\be \label{eq:riccati_kappa}
    - \alpha \kappa^2 \hat{F}_{\kappa} - \hat{F}_{\kappa} \alpha \kappa^2 - \hat{F}_{\kappa}^2 + \gamma^2 = 0.
\ee 
The unique solution $\hat{F}_{\kappa}>0$ to \eqref{eq:riccati_kappa} is $\hat{F}_{\kappa} = - \alpha \kappa^2 + \sqrt{\alpha^2 \kappa^4 + \gamma^2},$ and the optimal state feedback gain is thus given by $K^{\rm LQR}= \mathcal{F} M_{\hat{K}^{\rm LQR}} \mathcal{F}^{-1}$ where
\be \hat{K}^{\rm LQR}_{\kappa} = - \hat{F}_{\kappa} = \alpha \kappa^2 - \sqrt{\alpha^2 \kappa^4 + \gamma^2}.
\ee 

We next obtain the input-output mapping specified by the controller implementation \eqref{eq:opt_implementation} using \eqref{eq:recover_K} as 
\be 
    \hat{K}^{\mathcal{H}_2}_{\kappa}(s) = \frac{\hat{\Phi}^{u}_{\kappa}(s)}{\hat{\Phi}^{\psi}_{\kappa}(s)} = \alpha \kappa^2
    - \sqrt{\gamma^2 + \alpha^2 \kappa^4}.
\ee 
This is indeed equal to the static state feedback gain $\hat{K}^{LQR}$, thus $K^{\rm LQR}= K^{\mathcal{H}_2}=K^{\rm opt}.$
\end{proof}

We know that $K^{\rm opt}$ will be a spatial convolution operator 
\be 
    u(t,\theta) = \int_{\xi \in \mathbb{T}} k(\theta- \xi) \psi(t, \xi) dt,
\ee 
where $k(\theta) = (\mathcal{F}^{-1} \hat{K})(\theta).$ Numerical integration shows that $k(\theta)$  exhibits spatial decay away from the origin (see Fig.~\ref{fig:kernel_comparisons}). 

 We compare the static controller $K^{\rm opt}$ to the dynamic controller implementation \eqref{eq:opt_implementation} depicted in Fig. \ref{fig:controller_implementation_diagram}. As previously noted, the transfer function $ \left( I - G(s) {\Phi}^{\psi}(s)\right)$ is strictly proper, thus its corresponding convolution kernel has only a dynamic part, which we denote as $f_{\Phi^\psi}(t, \theta)$. The transfer function $G(s) \Phi^u(s)$  has a corresponding convolution kernel with both a static and dynamic part, which we denote as $h_{\Phi^u}(\theta), f_{\Phi^u}(t, \theta)$, respectively. The static (direct feedthrough) component of $G(s)\Phi^u(s)$ is obtained by evaluating its Fourier symbol in the limit as $s \to \infty$, 
$$  \mathcal{F}^{-1}\lim_{s\to\infty} \hat G(s,\kappa)\hat \Phi^u(s,\kappa)
= h_{\Phi^u}(\theta)=(\mathcal{F}^{-1}\hat{K}^{\rm opt})(\theta), $$
which shows that the static part of the implementation \eqref{eq:opt_implementation} recovers the static LQR gain, as expected.  The ``dynamic components" of the systems $G(s) \Phi^u(s)$ and $ \left( I - G(s) {\Phi}^{\psi}(s)\right)$ cancel each other out, leaving only the static feedback policy $K^{\rm opt}.$ This is confirmed in the plot of $h_{\Phi^u}(\theta), f_{\Phi^u}(t, \theta),f_{\Phi^\psi}(t, \theta)$ shown in Fig.~\ref{fig:kernel_comparisons}.

\section{Incorporation of Spatio-Temporal Decay Rate Constraints} \label{sec:constraints}
Without constraints, the closed-loop formulation of the $\mathcal{H}_2$ problem \eqref{eq:CL_formulation} is equivalent to a standard LQR problem. 
The advantage of this closed-loop formulation is that it remains convex when spatio-temporal constraints are incorporated. A similar problem has been studied in \cite{bamieh2005convex}, where necessary and sufficient conditions for convexity were dependent on ``funnel causality" properties of the underlying system and desired constraint. Notably, our formulation does not rely on these conditions and thus admits a broader range of constraints while preserving convexity.

An added constraint that $ I- G\Phi^{\psi}, G\Phi^u \in \mathcal{S}$ enforces a structural constraint on the components of our controller. In this case, formulation \eqref{eq:CL_formulation} becomes: 
\be \begin{aligned} \label{eq:constrained}
 &\argmin_{\Phi^{\psi}, \Phi^u\in H_{\infty} }~ \left\| C \Phi^{\psi}  + D \Phi^u\right\|_{\mathcal{H}_2} \\
    &~~~~~{\rm s.t.}~~~~~~\Phi^{\psi}, \Phi^u \text{ spatially invariant}\\
    & ~~~~~~~~~~~~~~~ \Phi^{\psi}, \Phi^u \text{ satisfy \eqref{eq:affine_innerProduct}}, \\
    & ~~~~~~~~~~~~~~~  ~ I- G\Phi^{\psi}, G\Phi^u \in \mathcal{S} .
\end{aligned} \ee  
So long as the set $ \mathcal{S} $ defines a convex constraint on $\Phi^u, \Phi^{\psi}$ the constrained problem \eqref{eq:constrained} will be convex. 

Of particular interest are constraints that specify a degree of spatio-temporal locality on $\Phi^{\psi}$ and $\Phi^u$. In the finite-dimensional, discrete-time, finite-horizon setting, such constraints take the form of spatio-temporal FIR constraints on the closed-loop system \cite{wang2019system}. In the continuous spatial and temporal domain setting, we believe it is more natural to impose a level of decay, or ``approximate locality". This choice might be viewed as analogue to the ``leaky" SLS problem presented in \cite{matni2017scalable}. 

One possible spatio-temporal constraint is when $\mathcal{S}$ is the set of linear PDEs whose Green's functions $\rho(\cdot, \cdot)$ satisfy 
\be\label{eq:constraint}\|\rho(\cdot,t)-\rho_T(\cdot,t)\|_{L_2(\mathbb{T})}<\epsilon(t)\ee where $\rho_T$ is the spatial truncation of $\rho$ given by $$ 
\rho_T(\theta,\cdot)=\begin{cases}
\rho(\theta,\cdot ) &|\theta|<r\\
0 & {\rm otherwise}
\end{cases}$$ and $r\le\pi$ is the truncation radius that encodes communication limitations. This constraint ensures that the `tails' of $\rho$  are sufficiently small, potentially enforcing spatial locality. Choosing $\epsilon(t)$ to grow in $t$ allows $\rho$ to be less localized in space as time increases. While constraining $\rho$ to have compact support as in \cite{conger2024convex} admits desirable computational properties, the decay constraint \eqref{eq:constraint} is perhaps more physically relevant. 
This set $\mathcal{S}$ is convex, and guarantees that  a spatially truncated version of controller will perform well using small gain arguments. 

\color{black}

%In addition to providing optimal control policies for distributed parameter systems, we posit that solutions to these constrained optimization problems \eqref{eq:constrained} can provide insight to tradeoffs that occur in finite dimensional systems. In particular, if spatio-temporal decay rates are imposed on the closed-loop system's Green's function, the resulting controller will take the form of a spatio-temporal convolution operator whose kernel possesses these decay properties. A discretization of this spatial decay could be interpreted as a spatial ``locality structure" while the temporal decay could be interpreted as a ``memory limitation". Thus, a tradeoff between memory and locality requirements could be inferred. Indeed, such a tradeoff between locality and state dimension in distributed controllers is a relevant question in today's landscape that remains unanswered \cite{vamsi2015optimal, lessard2013structured}. 

Obtaining analytical and numerical solutions to constrained problems of the form \eqref{eq:constrained} is the subject of ongoing work. While analytical solutions will require clever formulations of the constraints, numerical solutions may leverage approaches such as those presented in \cite{conger2025convex} and benefit from computational approaches like \cite{shivakumar2020pietools}.

\section{Conclusion}
It was shown that the infinite-horizon LQR problem for the diffusion equation over the unit circle can be reformulated as an $\mathcal{H}_2$ problem over the system's closed-loop responses, analogous to the finite-dimensional System Level Synthesis methodology. An implementation of this optimal controller built from these closed-loop responses was shown to preserve internal stability.  It was further illustrated that this closed-loop formulation admits a class of spatio-temporal constraints on the closed-loop responses while preserving convexity. Generalizations of this work to general linear distributed parameter systems are the subject of ongoing work. Of particular interest include incorporation of boundary conditions through an embedding approach \cite{epperlein2016spatially} or alternate representations \cite{peet2021partial}; and comparison of analytical results from this work to the computational approach of \cite{conger2025convex}. 

\bibliography{references}
\bibliographystyle{ieeetr}

\begin{appendix}

    \subsection{Proof of Proposition~\ref{prop:projection}}
    \label{app:projection}
    \begin{proof}
    We replace the implicit constraint \eqref{eq:implicit_constraint} with an analogous explicit constraint using the following Lemma: 
\begin{lemma}\label{lem:explicit_constraint}
    Fix $\kappa \in \mathbb{Z}.$ The single-input single-output (SISO) transfer functions $\hat{\Phi}^{\psi}_{\kappa}, ~\hat{\Phi}^u_{\kappa}$ satisfy  $\hat{\Phi}^{\psi}_{\kappa}, ~\hat{\Phi}^u_{\kappa} \in \mathcal{RH}_{2}$  and \eqref{eq:implicit_constraint} 
    if and only if they are of the form 
    \be \begin{aligned} \label{eq:explicit_parameterization}
        \hat{\Phi}^{\psi}_{\kappa}(s) & = \frac{1}{s+1} \left(1 + \hat{\rho}_{\kappa}(s) \right),\\
        \hat{\Phi}^{u}_{\kappa}(s) & = \frac{\alpha \kappa^2 -1}{s+1} + \frac{s + \alpha \kappa^2}{s+1} \hat{\rho}_{\kappa}(s)
    \end{aligned} \ee 
    for some SISO transfer function $\hat{\rho}_{\kappa}(s) \in  \mathcal{RH}_{2}.$
\end{lemma}    

    \begin{proof}
         It is simple to verify sufficiency, so we only prove necessity. We write \eqref{eq:implicit_constraint} as \be \label{eq:phi_hat_u}\hat{\Phi}_\kappa^u(s)=\left((s+1)-(-\alpha\kappa^2+1)\right)\hat{\Phi}_\kappa^\psi(s)-1.
         \ee 
          Since $\hat{\Phi}_\kappa^\psi(s)\in \mathcal{RH}_2,$ it follows that $\hat{\Phi}_\kappa^u(s)\in \mathcal{RH}_2$ only if $\lim_{s\to\infty}(s+1)\hat{\Phi}_\kappa^\psi(s)=1.$ This means $\hat{\Phi}_\kappa^\psi(s)$ is of the form $\hat{\Phi}_\kappa^\psi(s)=\frac{1}{s+1}(1+\hat{\rho}_\kappa(s))$ for some $\hat{\rho}_\kappa(s)\in \mathcal{RH}_2.$ Plugging this into \eqref{eq:phi_hat_u} gives the explicit formulas \eqref{eq:explicit_parameterization}.
     \end{proof}

We can now use the parameterization \eqref{eq:explicit_parameterization} to write \eqref{eq:CLOpt_kappa} as an unconstrained $\mathcal{H}_2$ design problem over a single SISO transfer function decision variable $\hat{\rho}_{\kappa}:$
    \be \label{eq:opt_rho}
        \argmin_{\hat{\rho}_\kappa\in \mathcal{RH}_2}\Bigg\|\underbrace{\begin{bmatrix} \frac{\gamma}{s+1}\\\frac{\alpha\kappa^2-1}{s+1}\end{bmatrix}+\begin{bmatrix} \frac{\gamma}{s+1}\\\frac{\alpha\kappa^2+s}{s+1}\end{bmatrix}\hat{\rho}_\kappa}_{G(\hat{\rho}_\kappa):=H+N\hat{\rho}_\kappa}\Bigg\|^2_{\mathcal{H}_2}
    \ee 
where $H\in \mathcal{RH}_2, N=\in \mathcal{RH}_\infty$.
        We want to leverage the following result \cite{kucera2006h2}: 
        \begin{lemma} \label{lem:orthog}
            Let $E\in \mathcal{RH}_2,F\in \mathcal{RH}_\infty$ with $F$ inner and $F^\sim E\in \mathcal{RH}_2^\perp.$ Then $\forall T\in \mathcal{RH}_2,$ \[\|E-FT\|^2_{\mathcal{H}_2}=\|E\|^2_{\mathcal{H}_2}+\|T\|^2_{\mathcal{H}_2}.\]
        \end{lemma}
    \vspace{.1in}
    Let $N=U_iU_o$ be the inner-outer factorization of $N$ where \be U_i=\begin{bmatrix}
        \frac{\gamma}{s+\sqrt{\gamma^2+\alpha^2\kappa^4}}\\\frac{s+\alpha\kappa^2}{s+\sqrt{\gamma^2+\alpha^2\kappa^4}}
    \end{bmatrix},~U_0=\frac{s+\sqrt{\gamma^2+\alpha^2\kappa^4}}{s+1}\ee  so that $U_i,U_o\in \mathcal{RH}_\infty$ and 
    $G(\hat{\rho}_\kappa)=H+U_iU_o\hat{\rho}_\kappa.$ Then $U_i^\sim G(\hat{\rho}_\kappa)=U_i^\sim H+U_o\hat{\rho}_\kappa=
    Y+Y^\perp+U_o$ where $Y,Y^\perp$ are the projections of $U_i^\sim$ onto $\mathcal{RH}_2,\mathcal{RH}_2^\perp,$ respectively. This gives $U_o\hat{\rho}_\kappa=U_i^\sim G(\hat{\rho}_\kappa)-Y-Y^\perp$ so that
    \begin{align}
        G(\hat{\rho}_\kappa)&=\underbrace{(H-U_iY)}_{=:G_1}-U_i\underbrace{(Y^\perp -U_i^\sim G(\hat{\rho}_\kappa))}_{=:N_1}.
    \end{align}
    Note that $G_1\in \mathcal{RH}_2$, $N_1=-Y-U_o\hat{\rho}_\kappa\in \mathcal{RH}_2$, and $U_i^\sim G_1=U_i^\sim H-Y=Y^\perp\in \mathcal{RH}_2^\perp$ so that by Lemma~\ref{lem:orthog}, 
    \begin{align}
       & \argmin_{\hat{\rho}_\kappa\in \mathcal{RH}_2}\|G(\hat{\rho}_\kappa)\|^2_{\mathcal{H}_2}=\argmin_{\hat{\rho}_\kappa\in \mathcal{RH}_2}\|G_1-U_iN_1\|^2_{\mathcal{H}_2}\\&~~
    = \argmin_{\hat{\rho}_\kappa\in \mathcal{RH}_2}\|G_1\|^2_{\mathcal{H}_2}+\|N_1\|^2_{\mathcal{H}_2}
    = \argmin_{\hat{\rho}_\kappa\in \mathcal{RH}_2}\|N_1\|^2_{\mathcal{H}_2}
    \end{align}
    since $G_1$ does not depend on $\hat{\rho}_\kappa$. To solve we set $N_1=0$ resulting in $\hat{\rho}_{\kappa}^{\rm opt.} (s) =U_o^{-1}Y=-U_o^{-1}\cdot(U_i^\sim H)\Big|_{\mathcal{RH}_2}$ so that the solution to \eqref{eq:opt_rho} is
    \be 
        \hat{\rho}_{\kappa}^{\rm opt.} (s) =  \frac{1 - \sqrt{\gamma^2 + \alpha^2 \kappa^4}}{s + \sqrt{\gamma^2 + \alpha^2 \kappa^4}},
    \ee
    which we plug into \eqref{eq:explicit_parameterization} to see \eqref{eq:opt_Phi} solves \eqref{eq:CLOpt_kappa}.
    \end{proof}

    \subsection{Proof of Lemma~\ref{lem:KStabilizing_PhiStable}} 
    \label{app:Phi_in_Hinf}
 \begin{proof}   Spatial invariance of $K$ implies spatial invariance of $\Phi^{\psi}$ and $\Phi^u$ from the relationship \eqref{eq:Phi} and properties of spatially-invariant systems. Let $\Sigma (A^K, B^K, C^K, D^K)$ denote the system with operator-valued transfer function $K$. $\lim_{s \rightarrow \infty} K(s) = D^K$ so that by equation \eqref{eq:Phi}, $\lim_{s \rightarrow \infty}\Phi^{\psi}(s) = 0$ and $\lim_{s \rightarrow \infty}\Phi^{u}(s) = 0$. The systems with transfer functions $\Phi^{\psi}, \Phi^u$ and $\mathbb{F}(P;K)$ can be represented as linear, spatially-invariant distributed parameter systems $\Sigma(A_{\rm c.l.}, I , I, 0), $ $\Sigma(A_{\rm c.l.}, I, C_u, 0), $ and $\Sigma(A_{\rm c.l.}, I, C_{\rm c.l.}, 0),$ respectively. Since $K$ is internally stabilizing, $\{T_{A_{\rm c.l.}(t)}\}$ is an exponentially stable $C_0$-semigroup, and thus the transfer function representations of these systems - $\Phi^{\psi}, \Phi^u$ and $\mathbb{F}(P;K)$ - are well-defined on $\{s \in \mathbb{C};~ {\rm Re}(s) \ge 0\}.$  By Plancherel's Theorem, 
$\sup_{{\rm Re}(s) \ge 0} ~ \| \Phi^{\psi}(s) \|_{L^2(\mathbb{T}) \rightarrow L^2(\mathbb{T})}  = \sup_{{\rm Re}(s) \ge 0} ~ \| M_{\hat{\Phi}^{\psi}(s)} \|_{\ell^2(\mathbb{Z}) \rightarrow \ell^2(\mathbb{Z})},$ and similarly for $\Phi^u,$ so that the pointwise condition \eqref{eq:Hinf_frequency} is equivalent to \eqref{eq:Hinf_condition}. 

To show that \eqref{eq:Hinf_frequency} holds, note that exponential stability of $\{T_{A_{\rm c.l.}}(t)\}$ is equivalent to exponential stability of the $C_0$-semigroup generated by $M_{\hat{A}_{\rm c.l.}} = \mathcal{F} A_{\rm c.l.} \mathcal{F}^{-1}$. This is a semigroup of multiplication operators, whose symbols are of the form $e^{t \hat{A}_{\rm c.l.}(\kappa)},$ where $\hat{A}_{\rm c.l.}(\kappa)$ is a finite-dimensional matrix for each $\kappa.$ Thus, 
there exists $m, \eta >0$ for which $\|T_{M_{\hat{A}_{\rm c.l.}}}(t) \| \le m e^{-\eta t}.$ The operator norm of each multiplication operator in this semigroup is given by the $\infty$-norm of its symbol \cite{bamieh2002distributed}, so that for each $t \ge 0,$
\be \label{eq:exp_decay_bound}
    m e^{-\eta t} \ge \|T_{M_{\hat{A}_{\rm c.l.}}}(t) \| = \sup_{\kappa} \sigma_{\max} \{ e^{t \hat{A}_{\rm c.l.}(\kappa)}\},
\ee 
where $\sigma_{\max}(N)$ denotes the maximum singular value of a matrix $N.$ {In addition to internal stabilization for each fixed $\kappa$, uniform stabilizability \cite[Def. 4.9]{morris2020controller} of the family of systems $(\hat{A}(\kappa),\hat{B}(\kappa))$ as $|\kappa|\to\infty$ ensures that the supremum in \eqref{eq:exp_decay_bound} is finite. Since $\hat B(\kappa)\hat K(\kappa)$ is uniformly bounded and $\hat A(\kappa)=-\alpha\kappa^2 I$,  $\lim_{|\kappa|\to\infty} \hat A_{\rm c.l.}(\kappa)\to-\infty$, so $(\hat{A}(\kappa),\hat{B}(\kappa))$ is uniformly stabilizable.} For each fixed $\kappa,$ the $H_{\infty}$ norm of the finite-dimensional transfer matrix $\hat{\Phi}^{\psi}_{\kappa}(s) = (sI-\hat{A}_{\rm c.l.}(\kappa) )^{-1}$, given by
$
    \sup_{{\rm Re}(s) \ge 0} \sigma_{\rm max} \{ (sI-\hat{A}_{\rm c.l.}(\kappa) )^{-1} \}, 
$
is equal to the operator norm of $\mathcal{O}\in \mathcal{L}(L^2[0, \infty))$ defined by $(\mathcal{O}f)(t) = \int_0^t e^{(t- \tau) \hat{A}_{\rm c.l.}(\kappa)}f(\tau) d \tau.$ By \eqref{eq:exp_decay_bound}, 
\be \begin{aligned} \label{eq:app_computations} 
& \|(\mathcal{O}f)(t) \|  = \left \| \int_0^t e^{(t- \tau) \hat{A}_{\rm c.l.}(\kappa)}f(\tau) d \tau \right\| \\
& = \left \| e^{t \hat{A}_{\rm c.l.}(\kappa)}\int_0^t e^{- \tau\hat{A}_{\rm c.l.}(\kappa)}f(\tau) d \tau \right\| \\
& \le \|  e^{t \hat{A}_{\rm c.l.}(\kappa)} \| \cdot \sup_{0 \le \tau \le t} \| e^{- \tau\hat{A}_{\rm c.l.}(\kappa)} \| \cdot \|f \| \\
& = \|  e^{t \hat{A}_{\rm c.l.}(\kappa)} \| \cdot \|f \|~ \le m e^{-\eta t} \| f\|,
\end{aligned} \ee 
where each $\| \cdot \| $ in \eqref{eq:app_computations} is the $L^2[0, \infty)$ norm.
Then, $\|\mathcal{O}f \| \le \|f\| \int_0^{\infty} (m e^{-\eta t})^2 dt,$ so that the operator norm of $\mathcal{O}\in \mathcal{L}(L^2[0, \infty))$ is upper bounded by $\tfrac{m^2}{\eta}< \infty$, independent of $\kappa.$ Thus, $\sup_{\kappa \in \mathbb{Z}} \sup_{{\rm Re}(s) \ge 0} \hat{\Phi}^{\psi}_{\kappa}(s) \le \tfrac{m^2}{\eta}< \infty,$ verifying the first inequality in \eqref{eq:Hinf_frequency}. Next, by definition of $C, D$ and the relationship \eqref{eq:Phi_to_F}, we have that $\mathbb{F}(P;K)(\kappa, s) =  \begin{bmatrix} \gamma \Phi^{\psi}_{\kappa} (s)  \\ 0 \end{bmatrix} + \begin{bmatrix} 0 \\ 1 \end{bmatrix} \Phi^u_{\kappa}(s)$, so that 
\be \begin{aligned} \nonumber
    &\sup_{{\rm Re}(s) \ge 0} \sup_{\kappa \in \mathbb{Z}} \big| \hat{\Phi}^u_{\kappa}(s) \big| \le \left \| \mathbb{F}(P;K) -  \gamma \begin{bmatrix}  \Phi^{\psi}  \\ 0 \end{bmatrix} \right\|_{\mathcal{H}_2}\\
    & ~~~~~\le ~ \left \| \mathbb{F}(P;K) \right\|_{\mathcal{H}_2} + \gamma \sup_{{\rm Re}(s) \ge 0, \kappa \in \mathbb{Z}} \big| \hat{\Phi}^{\psi}_{\kappa}(s) \big| < \infty,
\end{aligned} \ee 
verifying the second inequality in \eqref{eq:Hinf_frequency}.\end{proof}

     \subsection{Proof of Lemma~\ref{lem:op_constraint}}
     \label{app:op_constraint}
    \begin{proof} First, we demonstrate the equivalence of \eqref{eq:affine_innerProduct} to \eqref{eq:affine_pointwise} for each $s \in \{s \in \mathbb{C};~ {\rm Re}(s) \ge 0\}$. By Plancherel's Theorem, for a fixed $s,$ \eqref{eq:affine_innerProduct} is equivalent to 
    $$ \left<\hat{\eta}, M_{\hat{\Omega}(s)} \hat{\vartheta} \right>_{\ell^2(\mathbb{Z})} - \left< \hat{\eta}, M_{\hat{\Phi}^u(s)} \hat{\vartheta} \right>_{\ell^2(\mathbb{Z})}
        = \left< \hat{\eta}, \hat{\vartheta} \right>_{\ell^2(\mathbb{Z})},$$
        for all $\hat{\eta}, \hat{\vartheta} \in \ell^2(\mathbb{Z}),$
        where $\hat{\Omega}_{\kappa}(s) = \left( (s+1) - (\hat{A}_{\kappa} + 1) \right) \hat{\Phi}^{\psi}_{\kappa}(s)$. Expanding these inner products over $\ell^2(\mathbb{Z})$ gives 
        $ 
           0 =  \sum_{\kappa \in \mathbb{Z}} \hat{\eta}_{\kappa}\left( \hat{\Omega}_{\kappa}(s) - \hat{\Phi}^u_{\kappa} - 1 \right) \hat{\vartheta}_{\kappa}.
        $
        As this sum holds for all $\hat{\eta}, \hat{\vartheta} \in \ell^2(\mathbb{Z}),$ we have that  
        $
            0 =  \hat{\Omega}_{\kappa}(s) - \hat{\Phi}^u_{\kappa} - 1, 
        $
     for all $\kappa \in \mathbb{Z}$ or equivalently $ \left((s+1) - (- \alpha \kappa^2 + 1) \right) \hat{\Phi}^{\psi}_{\kappa}(s) - \hat{\Phi}_{\kappa}^u(s) = 1$ for all $\kappa \in \mathbb{Z}$.  Next, we demonstrate that the pointwise condition \eqref{eq:affine_pointwise} holds. From \eqref{eq:Phi}, we have that  
    \be \begin{aligned}
       & \hat{\Phi}^{\psi}_{\kappa}(s) = \tfrac{1}{s - \hat{A}_{\kappa} - \hat{K}_{\kappa}(s)} = \tfrac{1}{(s+1) -(- \alpha \kappa^2 + 1) - \hat{K}_{\kappa}(s)},\\
       &\hat{\Phi}^{\psi}_{\kappa}(s) = \tfrac{\hat{K}_{\kappa}(s)}{s - \hat{A}_{\kappa} - \hat{K}_{\kappa}(s)} = \tfrac{\hat{K}_{\kappa}(s)}{(s+1) -(- \alpha \kappa^2 + 1) - \hat{K}_{\kappa}(s)},
   \end{aligned}  \ee 
   which we can verify satisfy \eqref{eq:affine_pointwise}. \end{proof}
     \subsection{Proof of Lemma~\ref{lem:Controller_Implementation}}
     \label{app:controller_implementation}
     \begin{proof}Consider the controller implementation defined in \eqref{eq:control_implementation} together with the plant  \eqref{eq:transfer_function_state}. 
     For fixed $\kappa$, the assumption $\lim_{s\to\infty}s\hat{G}_\kappa^{-1}(s)=1$ together with the affine identity \eqref{eq:affine_pointwise} implies
$\lim_{s\to\infty}\hat{G}_\kappa(s)\hat{\Phi}^\psi_\kappa(s)=1$,
hence $I-G(s)\Phi^\psi(s)$ is strictly proper, and the feedback loop in Figure~\ref{fig:controller_implementation_diagram} is well-posed.  The controller \eqref{eq:control_implementation} must  yield bounded internal states and outputs in response to any bounded additive exogenous disturbances (see Fig.~\ref{fig:controller_implementation_diagram}). To show this, fix $\kappa$ and examine the scalar closed-loop transfer functions from the disturbance inputs $\hat w_\kappa(s), \hat n_\kappa(s)$, as shown in Fig. \ref{fig:controller_implementation_diagram}, to the signals $\hat\psi_\kappa(s)$, $\hat u_\kappa(s)$, and $\hat v_\kappa(s)$: %From Fig. \ref{fig:controller_implementation_diagram}, 
     \begin{align}
         \hat v_\kappa(s) &= \hat \eta_\kappa(s) + \hat \psi_\kappa(s)  + \hat n_\kappa(s)\\
         \hat \eta_\kappa(s) &= (1 - \hat G_\kappa(s)\hat \Phi^\psi_\kappa(s))\hat v_\kappa(s)\\
         \hat u_\kappa(s) &= \hat G_\kappa(s)\hat \Phi^u_\kappa(s) \hat v_\kappa(s)\\
         \hat \psi_\kappa(s) &= \tfrac{1}{s} \big( -\alpha\kappa^2 \hat\psi_\kappa(s) + 
         \hat u_\kappa(s) + \hat w_\kappa(s)\big)
     \end{align}
     Then using the affine identity in \eqref{eq:affine_pointwise} together with routine calculations yields the closed loop transfer matrices from $(\hat n_\kappa(s), \hat w_\kappa(s))$ to $(\hat\psi_\kappa(s), \hat u_\kappa(s), \hat v_\kappa(s))$ specified in \eqref{eq:transfer_functions}.
    % \begin{align}
     %    \begin{bmatrix}
     %        \hat\psi_\kappa\\
     %        \hat u_\kappa\\
     %        \hat v_\kappa
     %    \end{bmatrix} = \begin{bmatrix}
     %        \hat\Phi^u_\kappa(s) & \hat\Phi^\psi_\kappa(s)\\
    %         (s + \alpha\kappa^2)\hat\Phi^u_\kappa(s) & \hat\Phi^u_\kappa(s)\\
    %         (s + \alpha\kappa^2)\hat G^{-1}_\kappa(s)& \hat G^{-1}_\kappa(s)
    %     \end{bmatrix} \begin{bmatrix}
    %         \hat n_\kappa\\
     %        \hat w_\kappa\\
    %     \end{bmatrix}.
  %   \end{align}
%\begin{align}
%    \frac{\hat \psi_\kappa (s)}{\hat{w}_\kappa (s)} &= \frac{\hat\Phi_\kappa^\psi(s)}{\hat\Phi_\kappa^\psi(s)(s +\alpha\kappa^2) - \hat\Phi_\kappa^u(s)} = \hat\Phi_\kappa^\psi(s),\\
%    \frac{\hat v_\kappa (s)}{\hat{w}_\kappa (s)} &= \frac{1}{\hat G_\kappa(s) \Big[\hat\Phi_\kappa^\psi(s)(s +\alpha\kappa^2) - \hat\Phi_\kappa^u(s)\Big]} = \frac{1}{\hat G_\kappa(s)},
%\end{align}
By spatial invariance, the corresponding operator-valued transfer functions, $\frac{{\bf V}(s)}{{\bf W}(s)}, \frac{\bPsi(s)}{{\bf W}(s)}$, etc. are multiplication operators in the spatial frequency domain, whose symbols are specified in \eqref{eq:transfer_functions}. The assumptions $G^{-1}(s)(sI - A) \in \mathcal{H}_{\infty}$, $G(s)$ spatially invariant, and $\Phi^\psi, \Phi^u \in \mathcal{H}_\infty$, spatially invariant and strictly proper, imply that each operator transfer function is in $\mathcal{H}_\infty$. %Therefore, all internal signals remain bounded under bounded inputs, establishing internal stability of the controller implementation.

The disturbance-to-control and disturbance-to-state transfer functions $\frac{\hat u_\kappa (s)}{\hat w_\kappa (s)} = \hat\Phi_\kappa^u(s)$, and $\frac{\hat \psi_\kappa (s)}{\hat{w}_\kappa (s)} = \hat\Phi_\kappa^\psi(s)$ respectively are both in $\mathcal{H}_\infty$ and strictly proper, and yield,
$
    \hat u_\kappa (s) =\frac{\hat\Phi_\kappa^u(s)}{\hat\Phi_\kappa^\psi(s)} \hat\psi_\kappa (s),
$
which recovers the frequency-domain controller implementation \eqref{eq:recover_K}. \end{proof}
     
    \end{appendix}

\end{document}